\documentclass[11pt]{article}
\usepackage{amssymb,amsmath,tikz,mathrsfs}
\usepackage[T1]{fontenc}
\def\Pic{{\it Pic}}
\def\Ker{{\it Ker}}
\def\Res{{\it Res}}
\def\ord{{\it ord}}
\def\Div{{\it Div}}
\def\Jac{{\it Jac}}
\def\div{{\it div}}

\newtheorem{lemma}{Lemma}
\newtheorem{theorem}{Theorem}
\newtheorem{proposition}{Proposition}
\newcounter{mysection}
\begin{document}
%\hfill Updated \today.

\vspace{5mm}
\centerline{\textsc{Inverse spectral problem for GK integrable systems.}}
\vspace{5mm}
\centerline{\it V.V.Fock}
\vspace{5mm}

Given a minimal bipartite graph on a torus the main construction of Goncharov and Kenyon \cite{GK} gives a map from an algebraic torus of dimension one less than the number of faces of the graph to the space of pairs (planar curve, a line bundle on it) called the action-angle map. The aim of this section is to solve the inverse problem, namely given a planar curve of genus $g$ and a line bundle on it of degree $g-1$ construct a point of the algebraic torus. 

The key observation is that the space of pairs is birationally isomorphic to the configuration space of a collection of complete flags in an infinite dimensional vector space invariant with respect to an action of a free Abelian group with two generators. Therefore the construction of coordinates in the space of configurations of flags introduced in \cite{FG} for flag configurations in finite dimensional spaces applies with minor modifications.

Before giving the formula we recall some basic facts about combinatorics of bipartite graphs. In the appendix we briefly sketch necessary background on planar algebraic curves and theta functions.

\refstepcounter{mysection}
\paragraph{\arabic{mysection}. Discrete Dirac operator and the action-angle map.\label{s:discrete}}
In this section we give a very concise introduction to the Goncharov-Kenyon construction. Let $\Gamma$ be a bipartite graph with equal number of black and white vertices embedded into a two-dimensional torus $T$. Denote the set of black (resp. white )vertices of $\Gamma$ by $B$ and $W$, and by $F$ the set of connected components of the complement to $\Gamma$ in $\Sigma$, called faces.

For any graph $\Gamma$ a \textit{discrete line bundle} on it is just an association of a one dimensional vector space $V_v$ to every vertex $v$ of $\Gamma$. A discrete connection on a line bundle is an association to every edge of the graph of an isomorphism between vector spaces corresponding to its ends. For a bipartite graph it amounts to a collection isomorphisms $A_e:V_{b(e)}\to V_{w(e)}$ for every edge $e$ between the vertices $b(e)$ and $w(e$). Having chosen a basis in the vector spaces every isomorphism becomes just a multiplication by a nonzero number which we also denote by $A_e$. A collection of numbers $\boldsymbol{A}=\{A_e\}$ is called a \textit{discrete connection form} and can be interpreted as a cocycle $\boldsymbol{A}\in Z^1(\Gamma,\mathbb{C})$. Changing bases in the vector spaces amounts to changing the cocycle by a coboundary and therefore the space of connections can be identified with the cohomology group $H^1(\Gamma,\mathbb{C}^\times)$. 

Since the graph $\Gamma$ is embedded into the torus $T$ we have an exact sequence
$$1\to  H^1(T,\mathbb{C}^\times) \to H^1(\Gamma,\mathbb{C}^\times)\stackrel{d}{\to} H^2(T/\Gamma,\mathbb{C}^\times)\to H^2(T,\mathbb{C}^\times)\to 1$$

The space $H^2(T/\Gamma,\mathbb{C}^\times)$ is just the space denoted by $\mathcal{X}_\Gamma$ of associations of nonzero complex numbers $\boldsymbol{x}=\{x_i|i\in F\}$ to faces of the graph $\Gamma$. The differential $d \boldsymbol{A}\in H^2(T/\Gamma,\mathbb{C}^\times)$ can be interpreted as a discrete curvature of the connection --- association to every face of composition of isomorphisms corresponding to its sides. Denote by $\mathcal{X}_\Gamma^1 \in \mathcal{X}_\Gamma$ the image of the map $d$. A points  $\boldsymbol{x}\in \mathcal{X}_\Gamma^1$ are collections of numbers on faces with product equal to one.

The exact sequence implies also that $H^1(\Gamma,\mathbb{C}^\times)$ is a principal $H^1(\Sigma,\mathbb{C}^\times)$-bundle over the base $\mathcal{X}_\Gamma^1$.

Denote by $\boldsymbol{R}\in H^2(T/\Gamma,\pm  1)$ a 2-cocycle associating $-1$ to every face with the number of sides divisible by 4 and $1$ otherwise. A \textit{Kasteleyn orientation} $\boldsymbol{K}=\{K_e|e\in E\}$ is a cochain such that $d \boldsymbol{K}=\boldsymbol{R}$. 

A \textit{Dirac operator} on the graph $\Gamma$ provided with a discrete connection $\boldsymbol{A}$ is a map $$\mathfrak{D}_{\boldsymbol{K}}[\boldsymbol{A}]\colon\oplus_{b}V_b\to\oplus_w V_w$$ from the sum of spaces associated to black vertices to the sum of spaces associated to white vertices. It is defined by its action of $V_b$ as 
$$\mathfrak{D}_{\boldsymbol{K}}[\boldsymbol{A}]\vert_{V_b}=\oplus_{e|b(e)=b}A_eK_e.$$

From now on we will consider graphs such that for a generic connection the Dirac operator is nondegenerate. In this case $\mathfrak{D}_{\boldsymbol{K}}[\boldsymbol{A}]$ degenerates on a subvariety of $H^1(\Gamma,\mathbb{C}^\times)$ of codimension one. 

Given a point $\boldsymbol{x}\in \mathcal{X}_
\Gamma^1$ the preimage of it by the differential $d$ is isomorphic (up to a shift) to the algebraic torus $H^1(T,\mathbb{C}^\times)$. The intersection of the degeneration locus of the Dirac operator $\mathfrak{D}_{\boldsymbol{K}}[\boldsymbol{A}]$ with this torus is an algebraic curve $\Sigma_0(\boldsymbol{x})$. This curve can be compactified (see appendix \ref{s:planar}) to a curve $\Sigma(\boldsymbol{x})$ called the \textit{spectral curve} corresponding to the point $\boldsymbol{x}$ of the phase space. The kernel of $\mathfrak{D}_{\boldsymbol{K}}[\boldsymbol{A}]$ extends to $\Sigma(\boldsymbol{x})$ by continuity and defines a line bundle dual to a bundle $\mathcal{L}(\boldsymbol{x})$ (of degree $g-1$ as it will be clear from the second part of this note). The map associating a pair ($\Sigma(\boldsymbol{x}),\mathcal{L}(\boldsymbol{x})$) to a point $\boldsymbol{x}\in \mathcal{X}^1_\Gamma$ is called the \textit{action-angle map}.

In coordinates the equation of the curve $\Sigma(\boldsymbol{x})$ can be written as $\det \mathfrak{D}_{\boldsymbol{K}}[\boldsymbol{\lambda}d^{-1}(\boldsymbol{x})]=0$, where $\boldsymbol{\lambda}=(\lambda,\mu)\in H^1(T,\mathbb{C}^\times)$ is a cohomology class of the torus $T$. It is a Laurent polynomial equation with a Newton polygon $\Delta_\Gamma$ independent on $\boldsymbol{x}$, see the section \ref{s:graphs} and \cite{GK}.

The aim of this paper is to construct the inverse map, i.e., to restore the point of $\mathcal{X}^1_\Gamma$ out of a plane algebraic curve and a line bundle on it.

\refstepcounter{mysection}
\paragraph{\arabic{mysection}. Bipartite graphs and Newton polygons.\label{s:graphs}}
Recall that a \textit{zig-zag loop} \cite{GK} on a bipartite graph $\Gamma \subset T$ is a closed path along edges of $\Gamma$ turning maximally right at every white vertex and maximally left at every black one. Denote the set of zig-zag loops by $Z$. Every zig-zag loop $\alpha\in Z$ represents a homology class $\boldsymbol{h}_{\alpha}\in H_1(T,\mathbb{Z})$. Since there are exactly two zig-zag loops with opposite orientation passing through every edge, the sum of all such classes vanishes and therefore they form sides of a unique convex polygon $\Delta_\Gamma\subset H^1(T,\mathbb{R})$ defined up to a shift.

Let the dual surface $\Sigma$ be a surface obtained by gluing disks to the graph $\Gamma$ along every zig-zag loop. By construction the graph $\Gamma$ is embedded into $\Sigma$. 

The boundary of every face $i\in F$ is a curve representing a cycle $\check{\boldsymbol{h}}_i\in H_1(\Sigma,\mathbb{Z})$. Define a skew-symmetric matrix with integral entries $$\varepsilon_{ij}= \langle \check{\boldsymbol{h}}_i,\check{\boldsymbol{h}}_j\check{\rangle},$$ with $i,j \in F$ and $\langle,\check{\rangle}$ the intersection pairing of cycles on $\Sigma$. This matrix can be considered as an exchange matrix defining the structure of a cluster seed on $\mathcal{X}_\Gamma$ and in particular a Poisson structure on it by $\{x_i,x_j\}=\varepsilon_{ij}x_ix_j$. Since $\sum_{i\in F}\check{\boldsymbol{h}}_j=0$, we have $\sum_i \varepsilon_{ij}=0$ for any $j$ and therefore the submanifold $\mathcal{X}^1_\Gamma\in \mathcal{X}_\Gamma$ is a Poisson submanifold. 

The exchange matrix $\varepsilon_{ij}$ can be defined combinatorially as the number of common edges of the faces $i$ and $j$ taken with signs. To be more precise denote by $e_1,\ldots,e_{l(i)}$ the sides of the face $i$ taken counterclockwise and starting from the edge going from a black to a white vertex. Then $\varepsilon_{ij}=\sum_{k|e_k \in \partial j }(-1)^k$.

\refstepcounter{mysection}
\paragraph{\arabic{mysection}. Discrete Abel map.\label{s:graphsAbel}}

Consider a free Abelian group $\mathbb{Z}^{Z}$ generated by the set of zig-zag loops $Z$. This group has a natural $\mathbb{Z}$-grading with all generators having degree one. 

Consider the universal cover $\tilde{T}$ of the torus $T$ and denote by $\tilde{\Gamma}$ the lift of the graph $\Gamma$ to $\tilde{T}$. Define a map $\boldsymbol{d}$ associating to every face and vertex of $\tilde{\Gamma}$ an element of the group $\mathbb{Z}^{Z}$. Fix a face $i_0$ and fix any any element of $\mathbb{Z}^{Z}$ of degree 0 to be the value of $\boldsymbol{d}(i_0)$. 

Once this choice is done to define value of $\boldsymbol{d}$ on another face $i\in F$ choose a path $\gamma$ on $T$ connecting an internal point of the face $i$ to an internal of the face $i_0$ and define $\boldsymbol{d}(i)$ by the formula
\begin{equation}\label{path}
\boldsymbol{d}(i)= \boldsymbol{d}(i_0)+\sum_{\alpha\in Z}\langle\gamma,\boldsymbol{h}_{\alpha}\rangle\alpha,
\end{equation}
where $\langle\gamma,\boldsymbol{h}_{\alpha}\rangle$ is the intersection index between the zig-zag loop $\alpha$ and the path $\gamma$. 

It is obvious that the value of $\boldsymbol{d}$ does not depend on the choice of the path $\gamma$ since the expression (\ref{path}) vanishes for any closed path.

To define the the value of the map $\boldsymbol{d}$ on vertices slightly deform the zig-zag loops in such a way that they intersect edges in the midpoints only (as shown on fig. \ref{fi:blackvertex}). The connected components of the complement to the deformed loops correspond to either faces or vertices. The value $\boldsymbol{d}(v)$ is defined by exactly the same formula as for $\boldsymbol{d}(i)$ with the path $\gamma$ going from the midpoint of the face $i_0$ to the vertex $v$.

It is clear from the construction that $\deg \boldsymbol{d}(i)=0$, $\deg \boldsymbol{d}(w)=-1$ and $\deg \boldsymbol{d}(b)=1$ for any face $i$ white vertex $w$ and black vertex $b$, respectively.

The  map $\boldsymbol{d}$ is not unique, but is defined up to a shift $\boldsymbol{d}(i)\mapsto \boldsymbol{d}(i)+\boldsymbol{c}$ (and the same for the vertices) with any $\boldsymbol{c}\in \mathbb{Z}^Z$ of degree 0.

One can equivalently define the map $\boldsymbol{d}$ by the following property. Consider an edge $e$ of $\Gamma$ connecting a black vertex $b$ and the white vertex $w$. Denote by $\alpha^+$ zig-zag loop going along $e$ from $b$ to $w$ and $\alpha^-$ going along $e$ in the opposite direction. Denote also by $i^+$ and $i^-$ the faces of $\Gamma$ to the left and to the right form $e$ (viewed from $b$). Then we have
\begin{equation}\label{Dbalance}
\begin{array}{l}
\boldsymbol{d}(w)=\boldsymbol{d}(b)-\alpha^+-\alpha^-,\\
\boldsymbol{d}(i^+)=\boldsymbol{d}(b)-\alpha^+,\\
\boldsymbol{d}(i^-)=\boldsymbol{d}(b)-\alpha^-,
\end{array}
\end{equation}
These relations allow to unambiguously define the value of $\boldsymbol{d}$ starting from its value on any face of vertex.

Another property of the map $\boldsymbol{d}$ can be easily seen on the figure \ref{fi:blackvertex}B:
$$\sum_{i\in F}\varepsilon_{ij}\boldsymbol{d}(i)=0$$ 

Observe that any class of $H_1(T,\mathbb{Z})$ can be mapped to $\mathbb{Z}^Z$ by $\boldsymbol{h} \mapsto \sum_\alpha\langle\gamma,\boldsymbol{h}_{\alpha}\rangle\alpha$. We will denote a class and its image by the same letter.

Observe that the map $\boldsymbol{d}$ is equivariant with respect to the action or the group $H_1(T,\mathbb{Z})$, namely $\boldsymbol{d}(i+\boldsymbol{h}) = \boldsymbol{d}(i)+\boldsymbol{h}$ and similarly for $\boldsymbol{d}(w)$ or $\boldsymbol{d}(b)$. Therefore this map descends to the map from the faces and vertices of the graph $\Gamma$ to the quotient $\mathbb{Z}^Z/H_1(T,\mathbb{Z})$ which we will call the \textit{discrete Abel map} and denote by the same letter $\boldsymbol{d}$.

We call a bipartite graph $\Gamma$ \textit{minimal} if the number of faces is twice the area $S_{\Delta_\Gamma}$ of the polygon $\Delta_\Gamma$ and it is called \textit{simple} if all classes $\boldsymbol{h}_{\alpha}$ are non-divisible in $H_1(T,\mathbb{Z})$. As it is shown in \cite{GK} and \cite{FM} any graph can be transformed to make it minimal and simple without changing the Newton polygon.

\refstepcounter{mysection}
\paragraph{\arabic{mysection}. General solution.\label{s:general}}
The aim of this section and of the paper in general is to give an explicit inversion formula for the action-angle map. Namely, given a Laurent polynomial $P(\lambda, \mu)$ of two variables with a Newton polygon $\Delta$ defining an algebraic curve $\Sigma\subset H^2(T,\mathbb{C}^\times)$, a line bundle $\mathcal{L}\in \Pic^{g-1}(\Sigma)$ and a bipartite graph $\Gamma$ with the same Newton polygon describe a point $\boldsymbol{x}\in \mathcal{X}^1_\Gamma$ such that this point gives the pair $(\Sigma,\mathcal{L})$ as the value of the action-angle map.

The main observation permitting to solve this problem is that one can identify the set $Z$ of zig-zag loops with the set of points at infinity $\Sigma\backslash\Sigma_0$ in such a way that the corresponding homology classes $\boldsymbol{h}_{\alpha}\in H_1(T,\mathbb{Z})$ coincide. This map is not canonical, but is defined up to permutation of points at infinity corresponding to the same side of the Newton polygon. Once such isomorphism is chosen, we can extend this identification to a grading preserving map $\mathbb{Z}^Z\to \Div(\Sigma)$. The restriction of this map to $H_1(T,\mathbb{Z})$ takes value in principal divisors $\div(\Sigma)$ and thus defines a map $\mathbb{Z}^Z/H_1(T,\mathbb{Z}) \to \Pic(\Sigma)$ . In what follows we will denote both zig-zag loops and the corresponding points at infinity by the same letter.

Let $\mathcal{L} \in \Pic^{g-1}(\Sigma)$ be a line bundle on $\Sigma$. Denote by $H$ the space of meromorphic sections of the bundle $\mathcal{L}$ holomorphic on $\Sigma_0$ and denote by $F_\alpha^i=\{\psi\in H|\ord_\alpha \psi\geq i\}$ be the subspace of $H$ of section having zero of order at least $i$ at the point $\alpha \in \Sigma\backslash \Sigma_0$. The collection of the spaces $F_\alpha^i$ for a given $\alpha$ form a complete flag in the space $H$.  Observe that the group $H_1(T,\mathbb{Z})$ acts on $H$ by $\boldsymbol{h}: \psi \mapsto \langle \boldsymbol{h},\boldsymbol{\lambda}\rangle \psi$ and preserves the flags, namely 
\begin{equation}\label{period}
\boldsymbol{h} F^i_\alpha=F_\alpha^{i+\langle\boldsymbol{h},\boldsymbol{h}_\alpha\rangle}.
\end{equation}

Given $\boldsymbol{d}=\{d_\alpha\}\in \mathbb{Z}^Z$ denote by $F^{\boldsymbol{d}}=\cap F_\alpha^{d_\alpha}$ the intersection of such subspaces. For a generic $\mathcal{L}$ by the Riemann-Roch theorem the dimension of this intersection is given by $\dim F^{\boldsymbol{d}}=\max(0,-\deg \boldsymbol{d})$. Observe that $F^{\boldsymbol{d}}\subset F^{\boldsymbol{d}'}$ if $d_\alpha\geqslant d_\alpha'$ for any $i$.

Associate to every white vertex $w$ the one-dimensional space $V_w=F^{\boldsymbol{d}(w)}$. Let now $\alpha_1,\ldots,\alpha_l$ be the zig-zag loops passing through a black vertex $b$ and let $w_1,\ldots w_l$ be the corresponding white vertices as shown on Fig. \ref{fi:blackvertex}. From (\ref{Dbalance}) it follows that $\boldsymbol{d}(w_k)=\boldsymbol{d}(b)-\alpha_k-\alpha_{k+1}$. Consider now the kernel of the map $V_b=\Ker\left(\oplus_k F^{\boldsymbol{d}(w_k)}\to F^{\boldsymbol{d}(b)-\sum \alpha_k}\right)$. This map is well defined since $F^{\boldsymbol{d}(w_k)}\subset F^{\boldsymbol{d}(w)-\sum_k \alpha_k}$ and the kernel has dimension 1 if $\mathcal{L}$ is generic since $\dim F^{\boldsymbol{d}(b)-\sum \alpha_k}=l-1$. Obviously for any $k$ this map restricts to a collection of maps $\tilde{A}_k:V_b\to V_{w_k}$. Multiplying every such map by the value of the Kasteleyn connection form $\{K_e\}$ we get a connection form $\{A_e\}$ on the graph $\tilde{\Gamma}$. Taking into account the periodicity (\ref{period}) we see that the monodromy $x_i$ of the connection form $\{A_e\}$ around faces of $\tilde{\Gamma}$ is $H_1(T,\mathbb{Z})$-periodic and thus defines a point $\mathcal{X}^1_\Gamma$.

To give an explicit formula for $x_i$ in terms of $\theta$-functions, observe that for a white vertex $w$ of $\tilde{\Gamma}$ the space $V_w=F^{\boldsymbol{d}(w)}$ is the span of $\psi_{\boldsymbol{d}(w)}(z)=\theta_q(z-t+\boldsymbol{d}(w))E_{\boldsymbol{d}(w)}(z)$. To compute the maps $V_{b}\to V_{w_k}$ we need to find the coefficients $A_k$ of the identity
\begin{equation}\label{Diraccover}
\sum_k \tilde{A}_k\theta_q(z-t+\boldsymbol{d}(b)-\alpha_k-\alpha_{k+1})E_{\boldsymbol{d}(b)-\alpha_k-\alpha_{k+1}}(z)=0
\end{equation}
Comparing to the Fay's identity (\ref{Fay}) one gets

\begin{equation}\label{connection}
\tilde{A}_k = \frac{E(\alpha_l,\alpha_{l+1})}{\theta_q(t+\boldsymbol{d}(b)-\alpha_l)\theta_q(t+\boldsymbol{d}(b)-\alpha_{l+1})}
\end{equation}
or using the notations of (\ref{Dbalance}) and taking into account the Kasteleyn form

$$A_e=K_e\frac{E(\alpha^+,\alpha^-)}{\theta_q(t+\boldsymbol{d}(i^+))\theta_q(t+\boldsymbol{d}(i^-))}$$

Now computing the monodromy around a face $i$ with $2l$ sides enumerated counterclockwise starting from a black vertex as shown on Fig. \ref{fi:blackvertex}B and taking into account that $E(\alpha^+,\alpha^-)=\theta_{q'}(\boldsymbol{d}(i^+)-\boldsymbol{d}(i^-))/\sqrt{\phi(\alpha^+)\phi(\alpha^-)}$ one gets

\begin{equation}\label{xcoord}
x_i=(-1)^{l(i)/2+1}\prod_j \left(\frac{\theta_{q'}(\boldsymbol{d}(j)-\boldsymbol{d}(i))}{\theta_q(t+\boldsymbol{d}(j))}\right)^{\varepsilon_{ij}}.
\end{equation}

Now we are ready to formulate the main
\begin{theorem}
The formula (\ref{xcoord}) gives a point $\boldsymbol{x}\in \mathcal{X}_\Gamma^1$ such that $\Sigma(\boldsymbol{x})=\Sigma$ and $\mathcal{L}(\boldsymbol{x})=\mathcal{L}$.
\end{theorem}

To prove the theorem to give an explicit formulas for 
\begin{enumerate}
\item the connection $d^{-1}(\boldsymbol{x})$, 
\item the connection representing the class $\boldsymbol{\lambda}\in H^1(T,\mathbb{C^\times})$
\item a map $z\mapsto \boldsymbol{\lambda}(z)$ of the curve $\Sigma \to H^1(T,\mathbb{C^\times})$,
\item a map of the bundle $\mathcal{L^*}$ to the trivial bundle over $\Sigma$ with fiber $\mathbb{C}^W$.
\end{enumerate}
and verify that the Dirac operator $\mathfrak{D}[\boldsymbol{\lambda}(z)d^{-1}(\boldsymbol{x})]$ degenerates on the image of the bundle $\mathcal{L^*}$.

Since the connection given by the formula (\ref{connection}) is obviously $H_1(T,\mathbb{Z})$-periodic it defines a connection on the graph $\Gamma$ and not only on its cover $\tilde{\Gamma}$. 

The map $\Sigma\to H^2(T,\mathbb{C}^\times)$ can be defined by the formula $\langle\boldsymbol{\lambda}(z),\boldsymbol{h}\rangle=E_{\boldsymbol{h}}(z)$.

Fix now a lift of the vertices of $\Gamma$ to its cover $\tilde{\Gamma}$. Using such lift we can define the value of the map $\boldsymbol{d}$ for vertices of the graph $\Gamma$, but the relations (\ref{Dbalance}) will be satisfied only up to $H_1(T,\mathbb{Z})$. In particular for every edge $e$ of $\Gamma$ we can associate a homology class $\boldsymbol{h}_e=\boldsymbol{d}(b)-\boldsymbol{d}(w)-\alpha^+-\alpha^-$. For a given class $\boldsymbol{\lambda}\in H^1(T,\mathbb{C}^\times)$ we can also define a connection form$\{B_e\}$ by the formula $B_e = \langle \boldsymbol{\lambda},\boldsymbol{h}_e\rangle$. 
This connection represents the class $\boldsymbol{\lambda}$. Indeed, consider a closed path $\gamma$ passing consecutively through the edges $e_1,\ldots,e_{2l}$ starting from an edge going from a black to the white vertex. The monodromy along $\gamma$ is given by 
$$\prod_k \langle \boldsymbol{\lambda},\boldsymbol{h}_{e_k} \rangle^{(-1)^k}= \langle \boldsymbol{\lambda},\sum_k(-1)^k\boldsymbol{h}_{e_k} \rangle=\langle\boldsymbol{\lambda},\sum_{\alpha}\langle\gamma,\alpha\rangle\alpha\rangle=\langle\boldsymbol{\lambda},\gamma\rangle$$
and thus this connection form represents the class $\boldsymbol{\lambda}$. The second equality is true since after the substitution of expression for $\boldsymbol{h}_e$ all $\boldsymbol{d}(b)$ and $\boldsymbol{d}(e)$ cancel and $\alpha$'s enters with a sign given by the intersection index with the edge $e$.

The map of the bundle $\mathcal{L^*}$ to the trivial bundle $\mathbb{C}^W\times \Sigma$ is given by a collection of the sections $\psi_{\boldsymbol{d}(w)}$ of the dual bundle $\mathcal{L}$ for $w$ in the chosen vertices of $\tilde{\Gamma}$.

Now we need to put all these together and verify that for any black vertex $b$ of $\Gamma$ having neighbors $w_1,\ldots,w_l$ connected by edges $e_1,\ldots,e_l$ and with zig-zag loops $\alpha_1,\ldots,\alpha_k$ we have to verify that $\sum_kK_eA_{e_k}B_{e_k}\psi_{\boldsymbol{d}(w_k)}=0$.

Observe now that $B_{e_k}\psi_{\boldsymbol{d}(w_k)}=\psi_{\boldsymbol{d}(b)-\alpha_k-\alpha_{k+1}}$ and therefore the identity to be proven is a consequence of the relation (\ref{Diraccover}).

\begin{figure}[ht]
\begin{center}
\begin{tikzpicture}[scale=1.5]
\fill (0,0) circle (0.3);
\foreach \a in {-1,0,1,2,3}{
\draw[thick] (0,0)--(\a*70:2); \fill (\a*70:2) circle (0.3);\fill[color=white] (\a*70:2) circle (0.28);
\draw[<-<] (\a*70-15:2)..controls(\a*70+35:0.5)..(\a*70+85:2);
}
\fill[color=gray!0] (230:0.5)--(230:2.5)--(280:2.5)--(330:2.5)--(330:0.5)--cycle;
\draw (280:1) node {$\cdots$};
\draw (15:2.3) node {$\alpha_{k-1}$}; \draw (0:2) node {\scriptsize$k-1$};
\draw (85:2.2) node {$\alpha_{k}$};\draw (70:2) node {\scriptsize${k}$};
\draw (155:2.3) node {$\alpha_{k+1}$};\draw (140:2) node {\scriptsize${k+1}$};
\draw (225:2.2) node {$\alpha_{k+2}$};\draw (210:2) node {\scriptsize${k+2}$};
\draw (35:1.5) node {$i_k$};\draw (105:1.5) node {$i_{k+1}$};\draw (175:1.5) node {$i_{k+2}$};
\draw (270:2.5) node {A};
\end{tikzpicture}\qquad
\begin{tikzpicture}[scale=1.3]
\draw[thick] (-140:2)--(-70:2);
\foreach \a in {-1,0,1,2}{
\draw[thick] (\a*70:2)--(\a*70+70:2); \draw[thick] (\a*70:2)--(\a*70-5:2.5); \draw[thick] (\a*70:2)--(\a*70+5:2.5); \draw[thick,fill=black] (\a*70:2) circle (0.3);
}
\draw[>->] (-70+50:2.1)--(-70-50:2.1);
\draw[>->] (0-50:2.1)--(0+50:2.1);
\draw[>->] (70+50:2.1)--(70-50:2.1);
\draw[>->] (140-50:2.1)--(140+50:2.1);
\draw[thick,fill=white] (-70:2) circle (0.3);
\draw[thick,fill=white] (70:2) circle (0.3);
\fill[color=white] (160:2.3) arc (160:270:2.3)--(0,0)--cycle;
\draw (-70:2) node {\scriptsize $k-1$};
\draw (0:2) node {\textcolor{white}{\bf \scriptsize $k$}};
\draw (70:2) node {\scriptsize $k+1$};
\draw (140:2) node {\textcolor{white}{\bf \scriptsize $k+2$}};
\draw (-90:1.5) node {$\alpha_{k-1}$};
\draw (15:2.5) node {$\alpha_{k+1}$};
\draw (52:2.3) node {$\alpha_{k}$};
\draw (165:1.4) node {$\alpha_{k+2}$};
\draw (0,0) node {$i$};
\draw (35:2.2) node {$i_k$};\draw (105:2.2) node {$i_{k+1}$};\draw (-35:2.2) node {$i_{k-1}$};
\draw (270:3) node {B};
\end{tikzpicture}

%\begin{tikzpicture}[scale=1.5]
%\foreach \a in {-1,0,1,2,3}{
%\draw[thick] (0,0)--(\a*70:2); \fill (\a*70:2) circle (0.3);
%\draw[->] (\a*70-15:2)..controls(\a*70+35:0.5)..(\a*70+85:2);
%}
%\fill[color=gray!0] (230:0.5)--(230:2.5)--(280:2.5)--(330:2.5)--(330:0.5)--cycle;
%\draw (280:1) node {$\cdots$};
%\draw (-15:2.3) node {$\check{f}_{i-1}$}; \draw (0:2) node[color=gray!0] {\scriptsize$\boldsymbol{i-1}$};
%\draw (55:2.2) node {$\check{f}_{i}$};\draw (70:2) node[color=gray!0] {\scriptsize$\boldsymbol{i}$};
%\draw (125:2.2) node {$\check{f}_{i+1}$};\draw (140:2) node[color=gray!0] {\scriptsize$\boldsymbol{i+1}$};
%\draw (195:2.3) node {$\check{f}_{i+2}$};\draw (210:2) node[color=gray!0] {\scriptsize$\boldsymbol{i+2}$};
%\draw[thick,black,fill=gray!0] (0,0) circle (0.3);
%\end{tikzpicture}
\end{center}
\caption{\label{fi:blackvertex}}
\end{figure}

\begin{figure}
\begin{center}
\begin{tikzpicture}
\draw[white] (-0.5,1) node {\scriptsize $a\!-\!b$};
\node[white] at (1,-0.4) {$y(1\!+\!x^{-1})^{-1}$};
\node[white] at (1,2.4) {$u(1\!+\!x^{-1})^{-1}$};
\node at (1,1) {$x$};
\node at (1,-0.2) {$y$};
\node at (-0.2,1) {$z$};
\node at (1,2.2) {$u$};
\node at (2.2,1) {$v$};
\draw[thick] (0,2.3)--(0,2)--(-0.3,2);
\node at (-0.1,2.3) {.};\node at (-0.2,2.2) {.};\node at (-0.3,2.1) {.};
\draw[thick] (0,-0.3)--(0,0)--(-0.3,0);
\node at (-0.1,-0.3) {.};\node at (-0.2,-0.2) {.};\node at (-0.3,-0.1) {.};
\draw[thick] (2.3,0)--(2,0)--(2,-0.3);
\node at (2.3,-0.1) {.};\node at (2.2,-0.2) {.};\node at (2.1,-0.3) {.};
\draw[thick] (2,2.3)--(2,2)--(2.3,2);
\node at (2.1,2.3) {.};\node at (2.2,2.2) {.};\node at (2.3,2.1) {.};
\draw[very thick] (0,2)--(2,2)--(2,0)--(0,0)--cycle;
\draw[thick,fill] (0,2) circle (0.08);
\draw[thick,fill] (2,0) circle (0.08);
\draw[thick,fill=white] (2,2) circle (0.08);
\draw[thick,fill=white] (0,0) circle (0.08);
\draw[<->,thick] (4,1)--(4.8,1);
\node[white] at (1,-0.2) {$y(1\!+\!x)$};
\node at (-2,1) {A:};
\end{tikzpicture}\qquad
\begin{tikzpicture}
\node at (1,1) {$x^{-1}$};
\node at (1,-0.4) {$y(1\!+\!x^{-1})^{-1}$};
\node at (-0.8,1) {$z(1\!+\!x)$};
\node at (1,2.4) {$u(1\!+\!x^{-1})^{-1}$};
\node at (2.8,1) {$v(1\!+\!x)$};
\draw[thick] (0,2.3)--(0,2)--(-0.3,2);
\node at (-0.1,2.3) {.};\node at (-0.2,2.2) {.};\node at (-0.3,2.1) {.};
\draw[thick] (0,-0.3)--(0,0)--(-0.3,0);
\node at (-0.1,-0.3) {.};\node at (-0.2,-0.2) {.};\node at (-0.3,-0.1) {.};
\draw[thick] (2.3,0)--(2,0)--(2,-0.3);
\node at (2.3,-0.1) {.};\node at (2.2,-0.2) {.};\node at (2.1,-0.3) {.};
\draw[thick] (2,2.3)--(2,2)--(2.3,2);
\node at (2.1,2.3) {.};\node at (2.2,2.2) {.};\node at (2.3,2.1) {.};
\draw[very thick] (0.3,1.7)--(1.7,1.7)--(1.7,0.3)--(0.3,0.3)--cycle;
\draw[very thick] (0,0)--(0.3,0.3);
\draw[very thick] (0,2)--(0.3,1.7);
\draw[very thick] (2,0)--(1.7,0.3);
\draw[very thick] (2,2)--(1.7,1.7);
\draw[thick,fill] (0,2) circle (0.08);
\draw[thick,fill] (2,0) circle (0.08);
\draw[thick,fill=white] (2,2) circle (0.08);
\draw[thick,fill=white] (0,0) circle (0.08);
\draw[thick,fill=white] (1.7,0.3) circle (0.08);
\draw[thick,fill=white] (0.3,1.7) circle (0.08);
\draw[thick,fill] (0.3,0.3) circle (0.08);
\draw[thick,fill] (1.7,1.7) circle (0.08);
\end{tikzpicture}
\end{center}
\begin{center}
\begin{tikzpicture}
\draw[<-<] (-0.2,1.8)..controls(0.2,0.2)..(1.8,-0.2);
\draw[<-<] (2.2,1.8)..controls(1.8,0.2)..(0.2,-0.2);
\draw[<-<] (2.2,0.2)..controls(1.8,1.8)..(0.2,2.2);
\draw[<-<] (-0.2,0.2)..controls(0.2,1.8)..(1.8,2.2);
\draw (0.12,2.2) node {\scriptsize $a$};
\draw (1.88,2.23) node {\scriptsize $d$};
\draw (0.12,-0.2) node {\scriptsize $b$};
\draw (1.88,-0.2) node {\scriptsize $c$};
\draw (-0.5,1) node {\scriptsize $a\!-\!b$};
\draw (2.5,1) node {\scriptsize $c\!-\!d$};
\draw (1,2.3) node {\scriptsize$c\!-\!b$};
\draw (1,-0.3) node {\scriptsize$a\!-\!d$};
\draw (1,1) node {\scriptsize$c\!-\!b\!+\!a\!-\!d$};
%\draw[white] (0,0)--(2,2) node[midway,sloped,black] {$a\!+\!c$};
\draw[thick] (0,2.3)--(0,2)--(-0.3,2);
\node at (-0.1,2.3) {.};\node at (-0.2,2.2) {.};\node at (-0.3,2.1) {.};
\draw[thick] (0,-0.3)--(0,0)--(-0.3,0);
\node at (-0.1,-0.3) {.};\node at (-0.2,-0.2) {.};\node at (-0.3,-0.1) {.};
\draw[thick] (2.3,0)--(2,0)--(2,-0.3);
\node at (2.3,-0.1) {.};\node at (2.2,-0.2) {.};\node at (2.1,-0.3) {.};
\draw[thick] (2,2.3)--(2,2)--(2.3,2);
\node at (2.1,2.3) {.};\node at (2.2,2.2) {.};\node at (2.3,2.1) {.};
\draw[very thick] (0,2)--(2,2)--(2,0)--(0,0)--cycle;
\draw[thick,fill] (0,2) circle (0.08);
\draw[thick,fill] (2,0) circle (0.08);
\draw[thick,fill=white] (2,2) circle (0.08);
\draw[thick,fill=white] (0,0) circle (0.08);
\draw[<->,thick] (4,1)--(4.8,1);
\node at (-2,1) {B:};
\end{tikzpicture}\qquad
\begin{tikzpicture}
\draw (0.12,2.3) node {\scriptsize $a$};
\draw (1.88,2.33) node {\scriptsize $d$};
\draw (0.12,-0.3) node {\scriptsize $b$};
\draw (1.88,-0.3) node {\scriptsize $c$};
\draw (-0.2,1) node {\scriptsize $a\!-\!b$};
\draw (2.2,1) node {\scriptsize $c\!-\!d$};
\draw (1,2) node {\scriptsize$c\!-\!b$};
\draw (1,0) node {\scriptsize$a\!-\!d$};
\draw (1,1) node {\scriptsize $0$};
%\draw[white] (0,2)--(2,0) node[midway,sloped,black] {$b\!+\!d$};
\draw[<-<] (-0.2,1.9)..controls(1.67,1.67)..(1.9,-0.2);
\draw[<-<] (2.2,1.9)..controls(0.33,1.67)..(0.1,-0.2);
\draw[<-<] (-0.2,0.1)..controls(1.67,0.33)..(1.9,2.2);
\draw[<-<] (2.2,0.1)..controls(0.33,0.33)..(0.1,2.2);
\draw[very thick] (0,2.3)--(0,2)--(-0.3,2);
\node at (-0.1,2.3) {.};\node at (-0.2,2.2) {.};\node at (-0.3,2.1) {.};
\draw[very thick] (0,-0.3)--(0,0)--(-0.3,0);
\node at (-0.1,-0.3) {.};\node at (-0.2,-0.2) {.};\node at (-0.3,-0.1) {.};
\draw[very thick] (2.3,0)--(2,0)--(2,-0.3);
\node at (2.3,-0.1) {.};\node at (2.2,-0.2) {.};\node at (2.1,-0.3) {.};
\draw[very thick] (2,2.3)--(2,2)--(2.3,2);
\node at (2.1,2.3) {.};\node at (2.2,2.2) {.};\node at (2.3,2.1) {.};
\draw[very thick] (0.3,1.7)--(1.7,1.7)--(1.7,0.3)--(0.3,0.3)--cycle;
\draw[very thick] (0,0)--(0.3,0.3);
\draw[very thick] (0,2)--(0.3,1.7);
\draw[very thick] (2,0)--(1.7,0.3);
\draw[very thick] (2,2)--(1.7,1.7);

\draw[thick,fill=white] (0,2) circle (0.08);
\draw[thick,fill=white] (2,0) circle (0.08);
\draw[thick,fill=white] (2,2) circle (0.08);
\draw[thick,fill=white] (0,0) circle (0.08);
\draw[thick,fill] (0,2) circle (0.08);
\draw[thick,fill] (2,0) circle (0.08);
\draw[thick,fill=white] (2,2) circle (0.08);
\draw[thick,fill=white] (0,0) circle (0.08);
\draw[thick,fill=white] (1.7,0.3) circle (0.08);
\draw[thick,fill=white] (0.3,1.7) circle (0.08);
\draw[thick,fill] (0.3,0.3) circle (0.08);
\draw[thick,fill] (1.7,1.7) circle (0.08);
\node[white] at (2.7,1.2) {$v(1\!+\!x)$};
\node[white] at (-0.7,0.8) {$z(1\!+\!x)$};
\end{tikzpicture}
\end{center}
\caption{\label{fi:spider}Spider move. A: action on coordinates on $\mathcal{X}_\Gamma$. B: Action on the discrete Abel map $\boldsymbol{d}$.}
\end{figure}
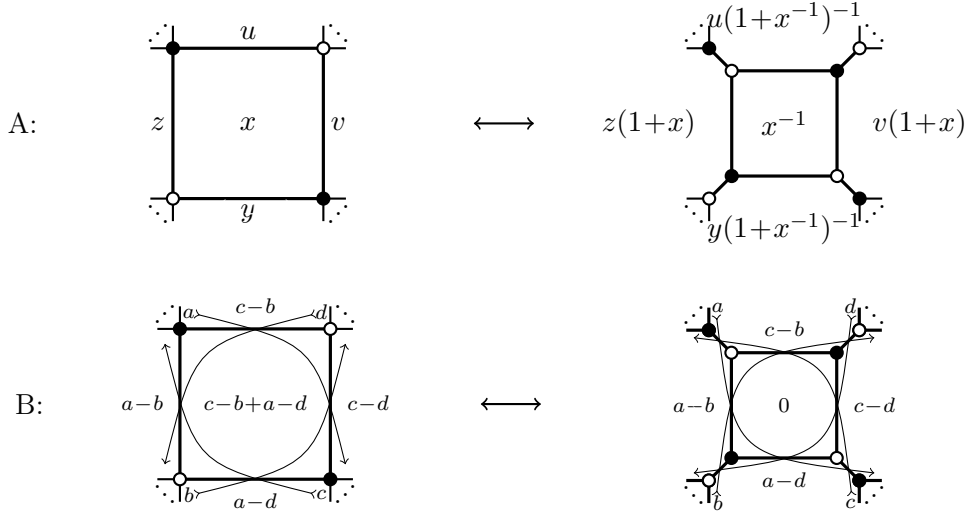

\refstepcounter{mysection}
\paragraph{\arabic{mysection}. Mutations.\label{s:mutations}}
As observed in \cite{GK} there exist transformations of the graph $\Gamma$ into another graph $\Gamma'$ that do not change the integrable system. Transformations of the first kind are retractions of two edges incident to a two-valent vertex or an inverse operation. The transformation of the second kind are applicable if there exists a quadrilateral face $i$ and it is shown on figure \ref{fi:spider}A. Transformation of the graphs induces birational isomorphisms between the corresponding tori $\mathcal{X}_\Gamma\to\mathcal{X}_{\Gamma'}$. For the transformation of the first kind this map is identity in coordinates (provided we use a natural bijection between faces of $\Gamma$ and $\Gamma'$). For the transformation of the second kind the isomorphism is a \textit{cluster mutation} and is shown also on fig. \ref{fi:spider}A. Define a correspondence between the discrete Abel maps before and after the transformation by the following conditions: For the transformation of the first kind we require its values on the corresponding faces coincide. For the transformation of the second kind its values coincide on corresponding faces except for the face $i$ and is changed to on $i$ to $\boldsymbol{d}(i)-a+b-c+d$, where $a,b,c,d$ are the zig-zag loops surrounding the quadrilateral face $i$ as is shown on fig. \ref{fi:spider}B. (In the figure we omitted $+\boldsymbol{d}(i')$ after every expression). 

Here we shall verify that the move is compatible with the parameterization given by the formula (\ref{xcoord}). More precisely we want to prove the following proposition.

\begin{proposition} 
The diagram
$$\begin{array}{rcl} 
&\Pic^{g-1}(\Sigma)&\\
\!\!\!\!&&\!\!\!\!\\
&\swarrow\qquad\qquad\searrow&\\[5pt]
\mathcal{X}_\Gamma& \longrightarrow &\mathcal{X}_{\Gamma'}
\end{array}
$$
is commutative for graph transformation of both kinds.
\end{proposition}

The proposition is obvious for the transformations of the first kind. For a spider move we need to verify that for any face of the graph $\Gamma'$ the coordinates given by (\ref{xcoord}) are related to the ones for the graph $\Gamma$ as shown on fig. \ref{fi:spider}A. 

$$x'=-\frac{E(a,b)}{\theta_q(b+c-t)}\frac{\theta_q(d+c-t)}{E(d,a)}\frac{E(c,d)}{\theta_q(d+a-t)}\frac{\theta_q(b+a-t)}{E(b,c)}=-\frac{F_t(a,b)F_t(c,d)}{F_t(d,a)F_t(b,c)}=x^{-1}
$$

$$y'=y\frac{E(b,c)}{\theta_q(a+c-t)}\frac{\theta_q(b+c-t)}{E(c,a)}\frac{E(d,a)}{\theta_q(b+d-t)}\frac{\theta_q(d+a-t)}{E(b,d)}=y\frac{F_t(b,c)F(d,a)}{F_t(c,a)F(b,d)}=$$
$$=y\frac{F_t(b,c)F_t(d,a)}{F_t(b,c)F_t(d,a)+F_t(a,b)F_(d,c)}=y(1+x^{-1})^{-1}
$$

The remaining relations are proven analogously

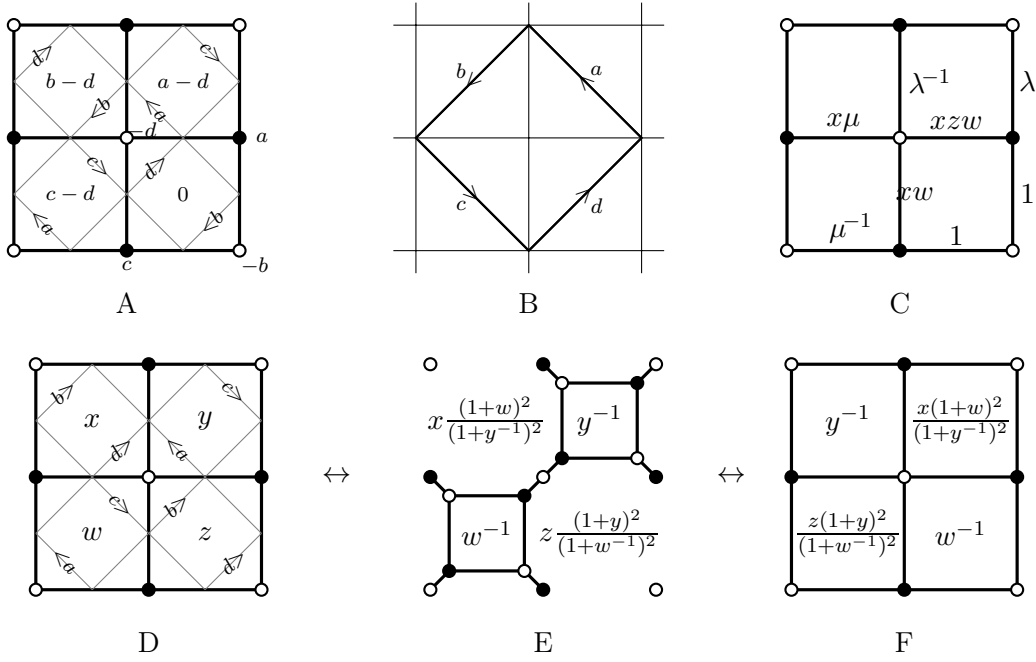
\begin{figure}
\begin{center}
\begin{tikzpicture}
\draw[very thick] (0,0)--(3,0)--(3,3)--(0,3)--cycle;
\draw[very thick] (0,1.5)--(3,1.5);
\draw[very thick] (1.5,0)--(1.5,3);
\draw[thick,fill=white] (0,0) circle (0.08);
\draw[thick,fill=white] (1.5,1.5) circle (0.08);
\draw[thick,fill=white] (3,3) circle (0.08);
\draw[thick,fill=white] (0,3) circle (0.08);
\draw[thick,fill=white] (3,0) circle (0.08);
\draw[thick,fill] (1.5,0) circle (0.08);
\draw[thick,fill] (0,1.5) circle (0.08);
\draw[thick,fill] (1.5,3) circle (0.08);
\draw[thick,fill] (3,1.5) circle (0.08);
\draw[thin,gray] (0.75,3)--(3,0.75) node[midway,sloped,black] {\scriptsize $<\!\!a$};
\draw[thin,gray] (3,0.75)--(2.25,0) node[midway,sloped,black] {\scriptsize $<\!\!b$};
\draw[thin,gray] (2.25,0)--(0,2.25) node[midway,sloped,black] {\scriptsize $c\!\!>$};
\draw[thin,gray] (0,2.25)--(0.75,3) node[midway,sloped,black] {\scriptsize $d\!\!>$};
\draw[thin,gray] (2.25,3)--(0,0.75) node[midway,sloped,black] {\scriptsize $<\!\!b$};
\draw[thin,gray] (0,0.75)--(0.75,0) node[midway,sloped,black] {\scriptsize $<\!\!a$};
\draw[thin,gray] (0.75,0)--(3,2.25) node[midway,sloped,black] {\scriptsize $d\!\!>$};
\draw[thin,gray] (3,2.25)--(2.25,3) node[midway,sloped,black] {\scriptsize $c\!\!>$};
\node at (0.75,2.25) {\scriptsize $b-d$};
\node at (2.25,2.25) {\scriptsize $a-d$};
\node at (2.25,0.75) {\scriptsize $0$};
\node at (0.75,0.75) {\scriptsize $c-d$};
\node at (1.7,1.6) {\scriptsize $-d$};
\node at (3.2,-0.2) {\scriptsize $-b$};
\node at (3.3,1.5) {\scriptsize $a$};
\node at (1.5,-0.2) {\scriptsize $c$};
\draw (1.5,-0.7) node {A};
\end{tikzpicture}\qquad\qquad
\begin{tikzpicture}
\draw[step=1.5,black,thin] (-1.8,-1.8) grid (1.8,1.8);
\draw[thick] (1.5,0)--(0,1.5)--(-1.5,0)--(0,-1.5)--cycle;
\draw[thick] (1.5,0)--(0,1.5) node[midway,sloped,black] {\scriptsize $<$};
\draw[thick] (0,1.5)--(-1.5,0) node[midway,sloped,black] {\scriptsize $<$};
\draw[thick] (-1.5,0)--(0,-1.5) node[midway,sloped,black] {\scriptsize $>$};
\draw[thick] (0,-1.5)--(1.5,0) node[midway,sloped,black] {\scriptsize $>$};
\draw (0.9,0.9) node {\scriptsize $a$}; 
\draw (0.9,-0.9) node {\scriptsize $d$};
\draw (-0.9,-0.9) node {\scriptsize $c$};
\draw (-0.9,0.9) node {\scriptsize $b$};
\draw (0,-2.2) node {B};
\end{tikzpicture}\qquad\qquad
\begin{tikzpicture}
\draw[very thick] (0,0)--(3,0)--(3,3)--(0,3)--cycle;
\draw[very thick] (0,1.5)--(3,1.5);
\draw[very thick] (1.5,0)--(1.5,3);
\draw[thick,fill=white] (0,0) circle (0.08);
\draw[thick,fill=white] (1.5,1.5) circle (0.08);
\draw[thick,fill=white] (3,3) circle (0.08);
\draw[thick,fill=white] (0,3) circle (0.08);
\draw[thick,fill=white] (3,0) circle (0.08);
\draw[thick,fill] (1.5,0) circle (0.08);
\draw[thick,fill] (0,1.5) circle (0.08);
\draw[thick,fill] (1.5,3) circle (0.08);
\draw[thick,fill] (3,1.5) circle (0.08);
\node at (3.2,2.25) {$\lambda$};
\node at (1.9,2.25) {$\lambda^{\!-1}$};
\node at (0.75,1.7) {$x\mu$};
\node at (0.85,0.3) {$\mu^{\!-1}$};
\node at (2.25,1.7) {$xzw$};
\node at (2.25,0.2) {$1$};
\node at (3.2,0.75) {$1$};
\node at (1.7,0.75) {$xw$};
\draw (1.5,-0.7) node {C};
\end{tikzpicture}
\end{center}
\begin{center}
\begin{tikzpicture}
\draw[very thick] (0,0)--(3,0)--(3,3)--(0,3)--cycle;
\draw[very thick] (0,1.5)--(3,1.5);
\draw[very thick] (1.5,0)--(1.5,3);
\draw[thick,fill=white] (0,0) circle (0.08);
\draw[thick,fill=white] (1.5,1.5) circle (0.08);
\draw[thick,fill=white] (3,3) circle (0.08);
\draw[thick,fill=white] (0,3) circle (0.08);
\draw[thick,fill=white] (3,0) circle (0.08);
\draw[thick,fill] (1.5,0) circle (0.08);
\draw[thick,fill] (0,1.5) circle (0.08);
\draw[thick,fill] (1.5,3) circle (0.08);
\draw[thick,fill] (3,1.5) circle (0.08);
\draw[thin,gray] (0.75,3)--(3,0.75) node[midway,sloped,black] {\scriptsize $<\!\!a$};
\draw[thin,gray] (3,0.75)--(2.25,0) node[midway,sloped,black] {\scriptsize $d\!\!>$};
\draw[thin,gray] (2.25,0)--(0,2.25) node[midway,sloped,black] {\scriptsize $c\!\!>$};
\draw[thin,gray] (0,2.25)--(0.75,3) node[midway,sloped,black] {\scriptsize $b\!\!>$};
\draw[thin,gray] (2.25,3)--(0,0.75) node[midway,sloped,black] {\scriptsize $d\!\!>$};
\draw[thin,gray] (0,0.75)--(0.75,0) node[midway,sloped,black] {\scriptsize $<\!\!a$};
\draw[thin,gray] (0.75,0)--(3,2.25) node[midway,sloped,black] {\scriptsize $b\!\!>$};
\draw[thin,gray] (3,2.25)--(2.25,3) node[midway,sloped,black] {\scriptsize $c\!\!>$};
%--(2.25,0)--(0,2.25)--cycle\draw[white] (0,2)--(2,0) node[midway,sloped,black] {$b\!+\!d$};
\node at (0.75,2.25) {$x$};
\node at (2.25,2.25) {$y$};
\node at (2.25,0.75) {$z$};
\node at (0.75,0.75) {$w$};
\node at (4,1.5) {$\leftrightarrow$};
\draw (1.5,-0.7) node {D};
\end{tikzpicture}\qquad
\begin{tikzpicture}
\draw[very thick] (0.25,0.25)--(1.25,0.25)--(1.25,1.25)--(0.25,1.25)--cycle;
\draw[very thick] (1.75,1.75)--(2.75,1.75)--(2.75,2.75)--(1.75,2.75)--cycle;
\draw[very thick] (0,0)--(0.25,0.25);
\draw[very thick] (1.25,1.25)--(1.75,1.75);
\draw[very thick] (2.75,2.75)--(3,3);
\draw[very thick] (0,1.5)--(0.25,1.25);
\draw[very thick] (1.5,0)--(1.25,0.25);
\draw[very thick] (1.5,3)--(1.75,2.75);
\draw[very thick] (3,1.5)--(2.75,1.75);
\draw[thick,fill] (0.25,0.25) circle (0.08);
\draw[thick,fill] (1.25,1.25) circle (0.08);
\draw[thick,fill] (1.75,1.75) circle (0.08);
\draw[thick,fill] (2.75,2.75) circle (0.08);
\draw[thick,fill=white] (0.25,1.25) circle (0.08);
\draw[thick,fill=white] (1.25,0.25) circle (0.08);
\draw[thick,fill=white] (1.75,2.75) circle (0.08);
\draw[thick,fill=white] (2.75,1.75) circle (0.08);
\draw[thick,fill=white] (0,0) circle (0.08);
\draw[thick,fill=white] (1.5,1.5) circle (0.08);
\draw[thick,fill=white] (3,3) circle (0.08);
\draw[thick,fill=white] (0,3) circle (0.08);
\draw[thick,fill=white] (3,0) circle (0.08);
\draw[thick,fill] (1.5,0) circle (0.08);
\draw[thick,fill] (0,1.5) circle (0.08);
\draw[thick,fill] (1.5,3) circle (0.08);
\draw[thick,fill] (3,1.5) circle (0.08);
\node at (0.75,2.25) {$x\frac{(1+w)^2}{(1+y^{-1})^2}$};
\node at (2.25,2.25) {$y^{-1}$};
\node at (2.25,0.75) {$z\frac{(1+y)^2}{(1+w^{-1})^2}$};
\node at (0.75,0.75) {$w^{-1}$};
\node at (4,1.5) {$\leftrightarrow$};
\draw (1.5,-0.7) node {E};
\end{tikzpicture}\quad
\begin{tikzpicture}
\draw[very thick] (0,0)--(3,0)--(3,3)--(0,3)--cycle;
\draw[very thick] (0,1.5)--(3,1.5);
\draw[very thick] (1.5,0)--(1.5,3);
\draw[thick,fill=white] (0,0) circle (0.08);
\draw[thick,fill=white] (1.5,1.5) circle (0.08);
\draw[thick,fill=white] (3,3) circle (0.08);
\draw[thick,fill=white] (0,3) circle (0.08);
\draw[thick,fill=white] (3,0) circle (0.08);
\draw[thick,fill] (1.5,0) circle (0.08);
\draw[thick,fill] (0,1.5) circle (0.08);
\draw[thick,fill] (1.5,3) circle (0.08);
\draw[thick,fill] (3,1.5) circle (0.08);
\node at (0.75,2.25) {$y^{-1}$};
\node at (2.25,2.25) {$\frac{x(1+w)^2}{(1+y^{-1})^2}$};
\node at (2.25,0.75) {$w^{-1}$};
\node at (0.75,0.75) {$\frac{z(1+y)^2}{(1+w^{-1})^2}$};
\draw (1.5,-0.7) node {F};
\end{tikzpicture}
\end{center}
\caption{A: Bipartite graph, zig-zag paths and value of the map $\boldsymbol{d}$ on faces and vertices; B: Newton polygon; C: Connection form $\{A_eB_e\}$; D:variables attached to faces; E:Mutaion of two faces; F: Genrerator of the discrete flow.\label{fi:Toda}}
\end{figure}

\refstepcounter{mysection}
\paragraph{\arabic{mysection}. Example.}
Consider the example of the simplest relativistic affine Toda lattice discussed in details in \cite{FM}. The bipartite graph $\Gamma$ is given on fig. \ref{fi:Toda}A, where we assume that the top of the picture is glued to the bottom and the left to the right. It has four zig-zag loops (shown on the fig. \ref{fi:Toda}C) representing the homology classes $(\pm 1,\pm 1)$ and defining the Newton polygon shown on fig. \ref{fi:Toda}B. We denote by $a,b,c,d$ the corresponding zig-zag loops. The vertical cycle of $H_1(T,\mathbb{Z})$ corresponds to $a+b-c-d$ and the variable $\lambda=E_{a+b-c-d}(z)$ while the horizontal one corresponds to $a+d-b-c$ and the variable $\mu=E_{a+d-b-c}(z)$. The labeling of vertices and faces by elements of $\mathbb{Z}^Z$ is shown on fig.\ref{fi:Toda}C. 

Since the curve $\Sigma$ is elliptic we have $\Pic^1(\Sigma)=\Sigma$ and $E(x,y)=\theta_{11}(x-y)$. Without loss of generality one can assume that $d=0$. The equations $a+d-b+c\in L_a$ and $a+b-c+d\in L_a$ has a solution $b=1/2$, $c=\tilde{p}_a+1/2$. Therefore
$$\lambda=\frac{\theta_{11}(z-a)\theta_{11}(z-1/2)}{\theta_{11}(z-a-1/2)\theta_{11}(z)};\quad \mu=\frac{\theta_{11}(z-a)\theta_{11}(z)}{\theta_{11}(z-a-1/2)\theta_{11}(z-1/2)}$$

The formula (\ref{xcoord}) gives

$$x=-\frac{\theta^2_{11}(1/2-a)}{\theta^2_{11}(a)}\frac{\theta^2_{00}(t+a)}{\theta^2_{00}(t+a+1/2)};\quad
y=-\frac{\theta^2_{11}(a)}{\theta^2_{11}(1/2-a)}\frac{\theta^2_{00}(t)}{\theta^2_{00}(t+1/2)}$$
$$z=-\frac{\theta^2_{11}(\tilde{p}_a)}{\theta^2_{11}(1/2-\tilde{p}_a)}\frac{\theta^2_{00}(t+1/2)}{\theta^2_{00}(t)};\quad
w=-\frac{\theta^2_{11}(1/2-a)}{\theta^2_{11}(a)}\frac{\theta^2_{00}(t+a+1/2)}{\theta^2_{00}(t+a)}$$

The group of discrete birational transformations is  $$\mathcal{G}_\Delta=\{(A,B,C,D)\in \mathbb Z^4|A+B+C+D=0\}/\langle (1,1,-1,-1),(1,-1,1,-1)\rangle = \mathbb{Z}\oplus (\mathbb{Z}/2\mathbb{Z})$$ 

One generator of this group is given by $$(x,y,z,w)\mapsto (z,w,x,y).$$ 
This generator has order two and corresponds to the automorphism of the graph.

The second generator of the infinite order corresponds to two mutations and shown on fig.\ref{fi:Toda}DEF:

$$(x,y,z,w)\mapsto \left(y^{-1}, x\frac{(1+w)^2}{(1+y^{-1})^2},w^{-1},z\frac{(1+y)^2}{(1+w^{-1})^2}\right)$$.

\newpage

\appendix
\noindent\textbf{\large{Appendix}}

\setcounter{mysection}{1}
\paragraph{\Alph{mysection}. Picard variety and the Abel map.\label{s:picard}}
Let $\Sigma$ be a smooth Riemann surface of genus $g$. The cohomology group $H^1(\Sigma,\mathbb{Z})$ possesses a nondegenerate skew-symmetric intersection form $\langle\cdot,\cdot\rangle$. Its complexification $H^1(\Sigma,\mathbb{C})$ contains a Lagrangian subspace $H^{1,0}(\Sigma,\mathbb{C})$ represented by holomorphic 1-forms on $\Sigma$. 

Fix two complimentary Lagrangian sublattices $L_a,L_b\subset H^1(\Sigma,\mathbb{Z})$ and denote their complexifications by $L_a^\mathbb{C}=L_a\otimes \mathbb{C}$ and $L^\mathbb{C}_b=L_b\otimes \mathbb{C}$, respectively. Three Lagrangian subspaces $L_a^\mathbb{C}$, $L_b^\mathbb{C}$ and $H^{1,0}(\Sigma,\mathbb{C})$ are transversal to each other. Denote by $\tau_a$ and $\tau_b$ the projections of $H^1(\Sigma,\mathbb{C})$ to $L_a^{\mathbb{C}}$ and $L_b^{\mathbb{C}}$, respectively, along $H^{1,0}(\Sigma,\mathbb{C})$. The skew-symmetric form $\langle\cdot,\cdot\rangle$ induces a bilinear form $Q$ on $H^1(\Sigma,\mathbb{C})$ by $Q(\boldsymbol{h}_1,\boldsymbol{h}_2)=\langle \tau_a( \boldsymbol{h}_1),\tau_b(\boldsymbol{h}_2)\rangle$. This form is symmetric since $\langle \tau_a(\boldsymbol{h}_1)-\tau_b(\boldsymbol{h}_1),\tau_a(\boldsymbol{h}_2)-\tau_b(\boldsymbol{h}_2)\rangle=0$. Its imaginary part is positive definite on $L_b^\mathbb{C}$ and negative definite on $L_a^\mathbb{C}$ since $0<i\langle\tau_a(\boldsymbol{h})-\tau_b(\boldsymbol{h}),\overline{\tau_a(\boldsymbol{h})-\tau_b(\boldsymbol{h}})\rangle=2\Im\langle\tau_a(\boldsymbol{h}),\overline{\tau_b(\boldsymbol{h})}\rangle$. The corresponding quadratic form is denoted by by the same letter $Q(\boldsymbol{h})=\langle \tau_a(\boldsymbol{h}), \tau_b(\boldsymbol{h})\rangle$. 

The quotient $L_a^{\mathbb{C}}/\tau_a (H^1(\Sigma,\mathbb{Z}))=L_b^{\mathbb{C}}/\tau_b (H^1(\Sigma,\mathbb{Z}))$ is called the Jacobian of the Riemann surface $\Sigma$ and denoted by $\Jac(\Sigma)$.

Let $\Div(\Sigma)$ be a free Abelian group generated by the set of points of the curve called the \textit{group of divisors} on $\Sigma$. It can also be considered as a group of singular $0$-chains on $\Sigma$. We make no distinction between a generator of the group $\Div(\Sigma)$ corresponding to a point $z\in\Sigma$ and the point $z$ itself and denote it by the same letter.  Any meromorphic section $\phi$ of a line bundle on $\Sigma$ induces a divisor $(\phi)=\sum_{z\in \Sigma} z\ord_z\phi$. The divisors of meromorphic functions are called \textit{principal}. The principal divisors form a subgroup of $\Div(\Sigma)$ denoted by $\div(\Sigma)$. The \textit{Picard variety} is the quotient $\Pic(\Sigma)=\Div(\Sigma)/\div(\Sigma)$. The group of divisors  has a grading $\Div(\Sigma)=\sum_i\Div^i(\Sigma)$ with any generator $z\in \Sigma$ having degree 1. This grading induces a grading $\Pic(\Sigma)=\sum_{i\in \mathbb{Z}}\Pic^i(\Sigma)$ since the grading of any principal divisor vanishes.

The variety $\Pic(\Sigma)$ can be also interpreted as the variety of line bundles on $\Sigma$. The correspondence associates to a given line bundle the divisor of any meromorphic section of it.

The \textit{Abel map} is an isomorphism $\mathcal{A}:\Pic^0(\Sigma)\to \Jac(\Sigma)$ defined as follows. For any $d\in \Pic^0(\Sigma)$ choose its representative $\tilde{d}\in\Div^0(\Sigma)$ and than choose a 1-chain $\partial^{-1}\tilde{d}$. These choices define a point $\mathcal{A}(d)$ in the space $L_b^{\mathbb{C}}$ by the property that $\langle \mathcal{A}(d),\omega\rangle=\int_{\partial^{-1}\tilde{d}}\omega$ for any holomorphic 1-form $\omega\in H^{1,0}(\Sigma)$. Making a different choice of $\partial^{-1}\tilde{d}$ or of $\tilde{d}$ results in changing the point $\mathcal{A}(d)$ by a point of the lattice $\tau_b(H^1(\Sigma,\mathbb{Z}))$ and thus $\mathcal{A}(d)$ considered as an element of $\Jac(\Sigma)$ is well defined.

In the text we make no distinction between $\Pic^0(\Sigma)$ and $\Jac(\Sigma)$. 

\refstepcounter{mysection}
\paragraph{\Alph{mysection}. Planar curves and Newton polygons.\label{s:planar}}
Let $P(\lambda,\mu)=\sum_{ij}c_{ij}\lambda^i\mu^j$ be a Laurent polynomial in two variables. Identify the pairs of numbers $(\lambda,\mu)=\boldsymbol{\lambda}$ with the cohomology classes from $H^1(T,\mathbb{C}^\times)$ of a two-dimensional torus $T$ with coefficients in the multiplicative group. The pairs of integers $(i,j)=\boldsymbol{h}$ can be than treated as homology classes from $H_1(T,\mathbb{Z})$ and monomials $\lambda^i\mu^j$ as natural pairing between cohomology and homology $\langle\boldsymbol{\lambda},\boldsymbol{h}\rangle$. 
The polynomial $P$ can be thus rewritten as $P(\boldsymbol{\lambda})=\sum_{\boldsymbol{\gamma}}c_{\boldsymbol{\gamma}}\langle \boldsymbol{\lambda},\boldsymbol{\gamma}\rangle$. Denote by $\Delta_P\subset H_1(\Sigma,\mathbb{R})$ the convex hull in of the set $\{\boldsymbol{\gamma}\in H_1(\Sigma,\mathbb{Z})|c_{\boldsymbol{\gamma}}\neq 0\}$. The equation $P(\boldsymbol{\lambda})=0$ defines an algebraic curve in $H^1(T,\mathbb{C}^\times)$.

Observe that the curve does not change if we make a transformation $P(\boldsymbol{\lambda})\to \langle\boldsymbol{\lambda},\boldsymbol{h}\rangle \beta P(\boldsymbol{\mu}\boldsymbol{\lambda})$, with any $\boldsymbol{h}\in H_1(T,\mathbb{Z}), \boldsymbol{\mu}\in H^1(\Sigma,\mathbb{C}^\times)$ and $\beta \in \mathbb{C}^\times$. Polynomials related by such transformations are called equivalent. The dimension of the space of equivalence classes of Laurent polynomials with given Newton polygon $\Delta_P$ is $I_{\Delta_P}+B_{\Delta_P}-3$, where $I_{\Delta_P}$ and $B_{\Delta_P}$ are the numbers of integer points strictly inside the polygon and on its boundary, respectively.

The curve is obviously non-compact but can be canonically compactified by adding points where either $\lambda$ or $\mu$ or both, considered as functions on the curve, vanish or have a pole.

The compactification goes as follows. Consider one side of the polygon $\Delta_P$. Without loss of generality we may assume that the side is a horizontal segment between $(0,0)$ and $(k,0)$ and $\Delta_P$ is located in the upper half plane. It means that the polynomial has the form $P(\lambda,\mu)=\sum_iP_i(\lambda)\mu^i$, where $P_i(\lambda)=0$ for $i<0$ and  $P_0(\lambda)$ is a polynomial of degree $k$ with nonvanishing degree zero term. The zero locus of $P(\lambda,\mu)$ intersects the line $\mu=0$ in the roots of the $P_0(\lambda)$. Adding the roots to the curve for every side of $\Delta_P$ makes it compact. The curve is nonsingular at the compactification point if the corresponding root is simple. We denote by $\Sigma$ the corresponding compact curve and call the added points \textit{the points at infinity}.

If the curve has no singularities at infinity the number of points corresponding to every side is equal to the number of segments into which a side is split by integer points (called \textit{boundary segments}). Remark that the bijection between points at infinity and boundary segments is not canonical and can be defined only up to permutations of the segments within every side. 

Consider a 1-form $\displaystyle \omega_{\boldsymbol{h}}=\Res\frac{\langle\boldsymbol{\lambda},\boldsymbol{h}\rangle}{P(\boldsymbol{\lambda})}\Omega$, where $\displaystyle \Omega=\frac{d\lambda\wedge d\mu}{\lambda\mu}$ be the canonical 2-form on $H^1(\Sigma,\mathbb{C}^\times)$ induced by the Poincar\'e pairing on the torus $T$. The form $\omega_{\boldsymbol{h}}$ is nonsingular and nonvanishing in the interior of $\Sigma$ provided the curve $\Sigma$ is nonsingular there. The zero order of the form $\omega_{\boldsymbol{h}}$ at a point $\alpha$ at infinity corresponding to the boundary segment $[\boldsymbol{x}_{\alpha},\boldsymbol{y}_{\alpha}]$ (with $\boldsymbol{x}_{\alpha}$ preceding $\boldsymbol{y}_{\alpha}$ in the counterclockwise order) is given by $\langle\boldsymbol{x}_{\alpha}-\boldsymbol{h},\boldsymbol{y}_{\alpha}-\boldsymbol{h}\rangle-1$, which is easy to verify assuming that the segment $[\boldsymbol{x}_{\alpha},\boldsymbol{y}_{\alpha}]$ belongs to the horizontal axis. Remark that the quantity $\langle\boldsymbol{x}_{\alpha}-\boldsymbol{h},\boldsymbol{y}_{\alpha}-\boldsymbol{h}\rangle$ is twice the signed area of the triangle $(\boldsymbol{x}_{\alpha},\boldsymbol{y}_{\alpha},\boldsymbol{h})$. It implies that the form $\omega_{\boldsymbol{h}}$ is holomorphic if $\boldsymbol{h}$ is strictly inside the polygon $\Delta_P$. It also implies that the degree of the canonical divisor $2g-2$ equals $2S_{\Delta_P}-B_{\Delta_P}$. By Pick's theorem $S_{\Delta_P}= I_{\Delta_P}+B_{\Delta_P}/2-1$ we get that $g=I_{\Delta_P}$ and thus the constructed holomorphic forms constitute a basis of $H^{1,0}(\Sigma,\mathbb{C})$.

\refstepcounter{mysection}
\paragraph{\Alph{mysection}. Parameterized curves and Newton polygons.\label{a:parameterized}}
 Let $\Sigma$ be an abstract smooth algebraic curve and let $\boldsymbol{\lambda}:\Sigma\to H^1(T,\mathbb{C}^\times)$ be a rational map to a two-dimensional algebraic torus which we consider as the cohomology group of a surface $T$ of genus 1. The points at infinity of this curve are those where the map is not defined. Denote the set of points at infinity by $Z$ enumerating them. For every $\alpha\in Z$ associate a homology class $\boldsymbol{h}_{\alpha}\in H_1(T,\mathbb{Z})$ by the condition that $\ord_{\alpha} \langle \boldsymbol{\lambda},\boldsymbol{h}\rangle\vert_{\Sigma}=\langle \boldsymbol{h},\boldsymbol{h}_{\alpha}\rangle$ for any $\boldsymbol{h}\in H_1(T,\mathbb{Z})$. Since the sum of orders of zeroes vanishes for any function we have $\sum_{\alpha\in Z} \boldsymbol{h}_{\alpha}=0$.  It implies that the vectors $\boldsymbol{h}_{\alpha}$ are sides of a convex polygon $\Delta_{\boldsymbol{\lambda}}\in H_1(T,\mathbb{R})$ defined up to a shift. 

We call the map $\boldsymbol{\lambda}$ nonsingular if its image has no singularities and if for any $\alpha\in Z$ the class $\boldsymbol{h}_{\alpha}$ is primitive (i.e., it is not divisible in $H^1(T,\mathbb{Z})$).

Now we are going to show that the image of the nonsingular map $\boldsymbol{\lambda}$ coincides with the curve given by the equation $P(\boldsymbol{\lambda})=0$ only if the corresponding Newton polygons $\Delta_{\boldsymbol{\lambda}}$ and $\Delta_P$ coincide.

Show first that every side of $\Delta_{\boldsymbol{\lambda}}$ is parallel to a side of $\Delta_P$. Indeed $\ord_{\alpha}P\circ\boldsymbol{\lambda}\geqslant\min_{\boldsymbol{h}\in\Delta_P}\langle \boldsymbol{h},\boldsymbol{h}_{\alpha}\rangle$ and the strict inequality can hold only if the minimum is attained at least twice. Since we need the order of $P\circ\boldsymbol{\lambda}$ at $\alpha$ to be $+\infty$ there exist at least two points $\boldsymbol{h}_1,\boldsymbol{h}_2\in \Delta_P$ such that $\langle \boldsymbol{h}_1-\boldsymbol{h}_2,\boldsymbol{h}_\alpha\rangle=0$ and $\langle \boldsymbol{h},\boldsymbol{h}_\alpha\rangle\geqslant\langle \boldsymbol{h}_1,\boldsymbol{h}_\alpha\rangle$ for any $\boldsymbol{h}\in \Delta$.

Since for a nonsingular map every point at infinity corresponds to a segment between integral points on the sides of $\Delta_{\boldsymbol{\lambda}}$, the lengths of the corresponding sides of $\Delta_P$ and $\Delta_{\boldsymbol{\lambda}}$ coincide and thus the polygons coincide also.

This observation allows to describe all maps $\Sigma\to H^1(T,\mathbb{C}^\times)$ with a given Newton polygon as collections of points $\alpha$ of $\Sigma$ one per segment $\boldsymbol{h}^\alpha$ between integer points of the boundary of $\Delta$ with the property 

$$\sum_{\alpha\in Z} \alpha\otimes \boldsymbol{h}_{\alpha}=0\in \Pic(\Sigma)\otimes H^1(\Sigma,\mathbb{Z}).$$

Since the dimension of the space of curves of genus $g$ is $3g-3$, the dimension of the space of pairs (curve, map) with a given Newton polytope is  $3g-3+B_{\Delta_{\boldsymbol{\lambda}}}-2g=g+B_{\Delta_{\boldsymbol{\lambda}}}-3$. On the other hand if $\boldsymbol{\lambda}$ is an isomorphism with the zero locus of $P$ this dimension must coincide with the dimension $I_{\Delta_P}+B_{\Delta_P}-3$ of equivalence classes of polynomials with a given Newton polytope $\Delta_P$. It gives another proof that $g=I_{\Delta_P}$.

Consider the subgroup  $\Div^\infty(\Sigma)\subset \Div(\Sigma)$ of divisors on $\Sigma$ supported at infinity. As an abstract group it is isomorphic to a free Abelian group generated by $Z$. The embedding $H_1(T,\mathbb{Z})\hookrightarrow \Div^\infty(\Sigma)$ given by $\boldsymbol{h}\mapsto (\langle\boldsymbol{\lambda},\boldsymbol{h}\rangle)$ induces an isomorphism between integer homology of $T$ and the set of principal divisors supported at infinity.

\refstepcounter{mysection}

\paragraph{\Alph{mysection}. Theta-functions.\label{a:theta}}

Recall that a spin structure on the curve $\Sigma$ can be identified with a quadratic form $q$ on $H^1(\Sigma,\mathbb{Z}/2\mathbb{Z})$ such that $q(\boldsymbol{h}^1+\boldsymbol{h}^2)=q(\boldsymbol{h}^1)+q(\boldsymbol{h}^2)+\langle \boldsymbol{h}^1,\boldsymbol{h}^2\rangle$ for any $\boldsymbol{h}^1,\boldsymbol{h}^2\in H^1(\Sigma,\mathbb{Z}/2\mathbb{Z})$. 

Given a decomposition $H^1(\Sigma,\mathbb{Z})=L_a\oplus L_b$ into two Lagrangian sublattices any such form can be presented as $q(\boldsymbol{l}_a+\boldsymbol{l}_b)=\langle \boldsymbol{l}_a,\boldsymbol{l}_b\rangle + \langle \boldsymbol{l}_a,\boldsymbol{\eta}\rangle + \langle \boldsymbol{\epsilon},\boldsymbol{l}_b\rangle$ for some $\boldsymbol{\eta} \in L_b/2L_b$ and $\boldsymbol{\epsilon} \in L_a/2L_a$.

If the surface $\Sigma$ has a complex structure spin structures can be identified with classes $q\in \Pic^{g-1}(\Sigma)$ such that $2q$ is the canonical class. The space $H^1(\Sigma,\mathbb{Z}/2\mathbb{Z})$ can be identified with classes $\boldsymbol{h}\in \Pic^0()\Sigma)$ such that $2\boldsymbol{h}=0$ and the class $q$ defines a quadratic form by $q(\boldsymbol{h})=\dim H^0(q+h)\mod{2}$ (see \cite{Atiyah}).

The theta function of characteristic $q$ is a function $L_a^{\mathbb{C}}\to \mathbb{C}$ defined by
$$ \theta_q(z)=\sum_{\boldsymbol{l}_b\in L_b}e^{2\pi i ({Q(\boldsymbol{l}_b+\boldsymbol{\eta}/2)}/{2}+\langle \boldsymbol{l}_b+\boldsymbol{\eta}/2,z+\boldsymbol{\epsilon}/2\rangle)} 
$$

The theta function satisfies the following straightforwardly verifiable properties:
\begin{enumerate}
\item $\theta_q(z+\boldsymbol{l}_a)=(-1)^{q(\boldsymbol{l}_a)}\theta_q(z)$ for any $\boldsymbol{l}_a\in L_a$,
\item $\theta_q(z+\tau_a(\boldsymbol{l}_b))=(-1)^{q( \boldsymbol{l}_b)}\theta_q(z)e^{-2\pi i ({Q(\boldsymbol{l}_b)}/{2}+\langle \boldsymbol{l}_b,z\rangle)}$ for any $\boldsymbol{l}_b\in L_b$,
\item $\theta_q(-z)=(-1)^{\langle \boldsymbol{\epsilon},\boldsymbol{\eta}\rangle}\theta_q(z)$.
\end{enumerate}

From these properties it follows that the zero locus of a theta function is a lift of a subvariety of $D[q]\in \Jac(\Sigma)$ called \textit{theta divisor}. It can be parameterized (it is the Jacobi inversion formula) by the $(n-1)$-st symmetric power $\Sigma^{n-1}/\mathfrak{S}_{n-1}$ by $(z_1,\ldots,z_{n-1})\mapsto q-\sum z_i$ or by $(w_1,\ldots,w_{n-1})\mapsto \sum w_i-q$.

The spin structure $q$ is called odd (respectively, even) if the theta function $\theta_q(z)$ is odd (respectively, even).

A map from $\Sigma$ to $\Pic^1(\Sigma)$ extends to a map of the Abelian universal cover $\widetilde{\Sigma}\to \widetilde{\Pic}^1(\Sigma)$. Let $e$ be any point of  $\widetilde{\Pic}^1(\Sigma)$. Consider the function on $\tilde{\Sigma}$ defined by $z\mapsto \theta_q(z-e)$. This function can be considered as a holomorphic section of a line bundle on $\Sigma$. Jacobi theorem implies that if $e=q-\sum_{i=1}^{g-2}z_i$ in $\Pic^1(\Sigma)$ then this section vanishes. Otherwise the divisor of this section is a sum of $g$ points and its class in $\Pic^g(\Sigma)$ is $q+e$.

\refstepcounter{mysection}
\paragraph{\Alph{mysection}. Prime form and Fay's trisecant identity.\label{s:prime}}
Recall following \cite{Fay} that the prime form is a function $E(x,y)$ on the square of the universal Abelian cover $\tilde{\Sigma}$ of the curve $\Sigma$ having the following properties:
\begin{enumerate}
\item $E(x,y)=-E(y,x)$
\item $E(x+\boldsymbol{l}_a,y)=E(x,
y)$ for any $\boldsymbol{l}_a\in L_a$
\item $E(x+\boldsymbol{l}_b,y)=e^{2\pi i (Q(\boldsymbol{l}_b)/2-\langle \boldsymbol{l}_b,x-y\rangle)}E(x,y)$ for any $\boldsymbol{l}_b\in L_b$
\item $E(x,y)=0$ if and only if $y=x+\boldsymbol{h}$ for some $\boldsymbol{h}\in H^1(\Sigma,\mathbb{Z})$.
\end{enumerate}
Here by $x+\boldsymbol{h}$ with $x\in \widetilde{\Sigma}$ and $\boldsymbol{h} \in H_1(\Sigma,\mathbb{Z})$ we mean the action of $\boldsymbol{h}$ on the point $x$.

Therefore the prime form as a function of its first argument $x$ can be considered as a section of the line bundle of degree 1 with divisor $y\in \Pic^1(\Sigma)$. 

The prime form allows to express explicitly a section of a line bundle on $\Sigma$ having a given divisor of zeroes and poles.Indeed, given a divisor $\boldsymbol{d}=\sum_\alpha d_\alpha \alpha \in \div(\Sigma)$ the product 
$$E_{\boldsymbol{d}}(z)=\prod_k E(z,\alpha)^{d_\alpha}$$ 
is a section of a line bundle corresponding to $\boldsymbol{d}$. In particular if $\boldsymbol{d}$ is principal than $E_{\boldsymbol{d}}$ is a single valued on $\Sigma$ with $(E_{\boldsymbol{d}})=\boldsymbol{d}$.

Let $q$ be an odd spin structure. According to the Abel theorem the divisor of the function $\theta_q(x-y)$, considered as a section of a line bundle on $\Sigma$ depending on $y\in \Sigma$ as a parameter, is equal to $y+\sum_{i=1}^{g-1}z_i$ with the class of $\sum_{i=1}^{g-1}z_i$ equal to $q$. 
Let $\phi$ be a holomorphic 1-form on $\Sigma$ with the divisor $2\sum_{i=1}^{g-1}z_i$. Then the prime form is given by 
\begin{equation}\label{prime}
E(x,y)=\theta_q(x-y)(\phi(x)\phi(y))^{-1/2}.
\end{equation}

The prime form does not depend on the choice of the spin structure $q$.

\begin{lemma}(Generalized Fay's trisecant identity)
Let $\{\alpha_k|k\in \mathbb{Z}/n\mathbb{Z}\}$ be a collection of points on a universal Abelian cover $\tilde{\Sigma}$ of a Riemann surface $\Sigma$ and $z$ be any other point of $\tilde{\Sigma}$. Let $t$ be any point of the universal cover of the Picard variety $\widetilde{\Pic}^1(\Sigma)$. Then the following identity holds: 
\begin{equation}\label{Fay}
\sum_k\frac{\theta(t+z-\alpha_k-\alpha_{k+1})}{E(z,\alpha_k)E(z,\alpha_{k+1})}\frac{E(\alpha_k,\alpha_{k+1})}{\theta(t-\alpha_k)\theta(t-\alpha_{k+1})}=0
\end{equation}

\end{lemma}

\noindent\textit{Proof:} Multiplying both sides by the common denominator one can reformulate the identity as
$$\sum_k \theta(t+z-\alpha_k-\alpha_{k+1})E(\alpha_k,\alpha_{k+1})\prod_{l\neq k,k+1}\theta(t-\alpha_l)E(z,\alpha_l)=0.
$$

Observe that all terms of the sum considered as functions of $z$ are holomorphic sections of a line bundle of degree $g+n-2$. Indeed, according to the Jacobi inversion formula the sum of zeroes in the Jacobian of every term (for generic $t$) is equal to $\sum \alpha_k-t$ and thus do not depend on $k$. We will show that the identity holds for $z=\alpha_m$. It implies that the sum of the remaining $g-2$ zeroes is equal to $-t$ for any $t\in \Pic^1(\Sigma)$, which is impossible for dimensional reasons.

To show that the identity holds for $z=\alpha_m$ we just
make this substitution and get
$$\sum_k \theta(t+\alpha_m-\alpha_k-\alpha_{k+1})E(\alpha_k,\alpha_{k+1})\prod_{l\neq k,k+1}\theta(t-\alpha_l)E(\alpha_m,\alpha_l)=$$
$$=\theta(t-\alpha_{m+1})E(\alpha_m,\alpha_{m+1})\prod_{l\neq m,m+1}\theta(t-\alpha_l)E(\alpha_m,\alpha_l)+$$
$$+\theta(t-\alpha_{m-1})E(\alpha_{m-1},\alpha_{m})\prod_{l\neq m-1,m}\theta(t-\alpha_l)E(\alpha_m,\alpha_l)
$$
$$=\prod_{l\neq m}\theta(t-\alpha_l)E(\alpha_m,\alpha_l)-\prod_{l\neq m}\theta(t-\alpha_l)E(\alpha_m,\alpha_l)=0,$$
where the first equality is satisfied since all but two terms of the sum vanish. Lemma is proven.

For $n=1$ and $2$ the lemma is trivial. For $n=3$ replacing $\alpha_0,\alpha_1,\alpha_2,z$ by $a,b,c,d$, respectively, and $t$ by $t+d$ we get the trisecant identity in a more usual form \cite{Mumford}:

$$\begin{array}{rl}
~&\theta(b+c-t)E(b,c)\theta(d+a-t)E(d,a)\ +\\
+&\theta(c+a-t)E(c,a)\theta(d+b-t)E(d,b)\ +\\
+&\theta(a+b-t)E(a,b)\theta(d+c-t)E(d,c)=0
\end{array}$$ 

Introducing a function $F_t(u,v)=\theta(u+v-t)E(u,v)$ this formula can be written in an even more elegant form
$$F_t(b,c)F_t(d,a)+F_t(c,a)F_t(d,b)+F_t(a,b)F_t(d,c)=0
$$

\end{document}